\documentclass[letterpaper]{IEEEtran}

\usepackage{mathdots}
\usepackage{amsthm,amsmath,amssymb}
\usepackage[hide links]{hyperref}
\usepackage[noabbrev]{cleveref}
\usepackage{cite}
\usepackage[inline]{enumitem}
\usepackage{tikz}
\usepackage{upgreek,stfloats,bm}
\usepackage[bb=boondox]{mathalfa}

\usepackage[no end]{algorithmic}
\usepackage{algorithm}

\usetikzlibrary{positioning,calc,patterns,decorations.pathmorphing,decorations.markings,shapes.geometric}

\newcommand{\D}{\mathrm{d}}
\newcommand{\R}{\mathbb{R}}
\newcommand{\Z}{\mathbb{Z}}
\newcommand{\mathe}{\mathrm{e}}
\newcommand{\C}{\mathbb{C}}
\newcommand{\N}{\mathbb{N}}

\newcommand{\Ker}{\mathrm{ker}}
\newcommand{\rank}{\mathrm{rank}}
\newcommand{\swap}{S}

\newcommand{\col}{\mathrm{col}}
\newcommand{\He}{\mathrm{He}}
\newcommand{\AEX}{A_{\mathrm{e}}}
\newcommand{\BEX}{B_{\mathrm{e}}}
\newcommand{\CEX}{C_{\mathrm{e}}}
\newcommand{\GEX}{G_{\mathrm{e}}}

\newcommand{\SP}{\mathrm{span}}

\newcommand{\TS}{T_S}

\newcommand{\rma}{\mathrm{a}}
\newcommand{\rmb}{\mathrm{b}}

\newcommand{\AREF}{A_{\mathrm{r}}}
\newcommand{\BREF}{B_{\mathrm{r}}}

\newcommand{\AIN}{A_{\mathrm{i}}}
\newcommand{\AINrr}{A_{\mathrm{i},\mathrm{r},\mathrm{r}}}
\newcommand{\AINru}{A_{\mathrm{i},\mathrm{r},\mathrm{u}}}
\newcommand{\AINuu}{A_{\mathrm{i},\mathrm{u},\mathrm{u}}}
\newcommand{\BIN}{B_{\mathrm{i}}}
\newcommand{\BINr}{B_{\mathrm{i},\mathrm{r}}}

\newcommand{\CIN}{C_{\mathrm{i}}}

\newcommand{\RIN}{R_{\mathrm{i}}}

\newcommand{\HF}{J_{\mathrm{f}}}

\newcommand{\AF}{A_{\mathrm{f}}}
\newcommand{\BF}{B_{\mathrm{f}}}

\newcommand{\GF}{G_{\mathrm{f}}}

\newcommand{\Au}{A_{\mathrm{r},m}}
\newcommand{\Ay}{A_{\mathrm{r},p}}
\newcommand{\Bu}{B_{\mathrm{r},m}}
\newcommand{\By}{B_{\mathrm{r},p}}
\newcommand{\AEF}{A_{\mathrm{M},\mathrm{f}}}
\newcommand{\BEF}{B_{\mathrm{M},\mathrm{f}}}

\newcommand{\GEF}{G_{\mathrm{M},\mathrm{f}}}

\newcommand{\bbe}[2]{\mathbb{e}_{#1}^{#2}}

\newcommand{\KF}{K_{\mathrm{f}}}

\newcommand{\AV}{A_{\mathrm{v}}}
\newcommand{\BV}{B_{\mathrm{v}}}
\newcommand{\GV}{G_{\mathrm{v}}}
\newcommand{\PIV}{\Pi_{\mathrm{v}}}

\newcommand{\PIF}{\Pi_{\mathrm{f}}}
\newcommand{\PIEX}{\Pi_{\mathrm{e}}}

\newcommand{\rmI}{\mathrm{I}}
\newcommand{\rmII}{\mathrm{II}}

\newcommand{\LF}{L_\mathrm{f}}
\newcommand{\Sk}{\mathrm{Sk}}

\newcommand{\AVD}{A_{\mathrm{v}}^{\mathrm{d}}}
\newcommand{\BVD}{B_{\mathrm{v}}^{\mathrm{d}}}

\newcommand{\epse}{\epsilon_{\mathrm{e}}}
\newcommand{\KEF}{K_{\mathrm{f}}}

\newcommand{\han}[1]{\mathfrak{H}_{#1}}

\newcommand{\AM}{A_{\mathrm{M}}}
\newcommand{\BM}{B_{\mathrm{M}}}
\newcommand{\CM}{C_{\mathrm{M}}}
\newcommand{\OM}{O_{\mathrm{M}}}
\newcommand{\MM}{M_{\mathrm{M}}}

\newcommand{\mD}{\mathcal{D}}

\newcommand{\AMIN}{A_{\mathrm{M},\mathrm{i}}}
\newcommand{\BMIN}{B_{\mathrm{M},\mathrm{i}}}
\newcommand{\CMIN}{C_{\mathrm{M},\mathrm{i}}}
\newcommand{\RMIN}{R_{\mathrm{M},\mathrm{i}}}
\newcommand{\AMINrr}{A_{\mathrm{M},\mathrm{i},\mathrm{r},\mathrm{r}}}
\newcommand{\AMINru}{A_{\mathrm{M},\mathrm{i},\mathrm{r},\mathrm{u}}}
\newcommand{\AMINuu}{A_{\mathrm{M},\mathrm{i},\mathrm{u},\mathrm{u}}}
\newcommand{\BMINr}{B_{\mathrm{M},\mathrm{i},\mathrm{r}}}

\newcommand{\AMEX}{A_{\mathrm{M},\mathrm{e}}}
\newcommand{\BMEX}{B_{\mathrm{M},\mathrm{e}}}
\newcommand{\CMEX}{C_{\mathrm{M},\mathrm{e}}}
\newcommand{\GMEX}{G_{\mathrm{M},\mathrm{e}}}

\newcommand{\PiM}{\Pi_{\mathrm{M}}}
\newcommand{\PIEF}{\Pi_{\mathrm{M},\mathrm{f}}}
\newcommand{\PIMEX}{\Pi_{\mathrm{M},\mathrm{e}}}
\newcommand{\AMV}{A_{\mathrm{M},\mathrm{v}}}
\newcommand{\BMV}{B_{\mathrm{M},\mathrm{v}}}
\newcommand{\GMV}{G_{\mathrm{M},\mathrm{v}}}
\newcommand{\AMVD}{A_{\mathrm{M},\mathrm{v}}^{\mathrm{d}}}
\newcommand{\BMVD}{B_{\mathrm{M},\mathrm{v}}^{\mathrm{d}}}
\newcommand{\PIMV}{\Pi_{\mathrm{M},\mathrm{v}}}

\theoremstyle{plain}
\newtheorem{theom}{Theorem}
\newtheorem{lemm}{Lemma}
\newtheorem{prpp}{Proposition}
\newtheorem{crl}{Corollary}
\newtheorem{asss}{Assumption}

\newtheorem{taggedassumptionx}{Assumption}
\newtheorem{taggedtheoremx}{Theorem}
\newtheorem{taggedlemmax}{Lemma}
\newtheorem{taggedpropositionx}{Proposition}

\newenvironment{taggedassumption}[1]
{\renewcommand\thetaggedassumptionx{#1}\taggedassumptionx}
{\endtaggedassumptionx}

\newenvironment{taggedtheorem}[1]
{\renewcommand\thetaggedtheoremx{#1}\taggedtheoremx}
{\endtaggedtheoremx}

\newenvironment{taggedlemma}[1]
{\renewcommand\thetaggedlemmax{#1}\taggedlemmax}
{\endtaggedlemmax}

\newenvironment{taggedproposition}[1]
{\renewcommand\thetaggedpropositionx{#1}\taggedpropositionx}
{\endtaggedpropositionx}

\theoremstyle{definition}

\theoremstyle{remark}
\newtheorem{remma}{Remark}

\crefname{prb}{problem}{problems}
\crefname{lemm}{lemma}{lemmas}
\crefname{asss}{assumption}{assumptions}
\crefname{theom}{theorem}{theorems}

\title{Derivative free data-driven stabilization \\of continuous-time linear systems\\ from input-output data}

\author{Corrado Possieri
\thanks{
	This work was supported in part by the Università degli Studi di Roma Tor Vergata through the Project FORM.\\ \indent
	This paper is an extended version of reference \cite{extended}. \\ \indent
	 The author is with the Dipartimento di Ingegneria Civile e Ingegneria Informatica, 
		Universit{\`a} degli Studi di Roma Tor Vergata, 00133 Roma, Italy (e-mail: 
		\href{mailto:corrado.possieri@uniroma2.it}{corrado.possieri@uniroma2.it}).}}

\begin{document}

\maketitle
	\begin{abstract}
This letter presents a data-driven framework for the design of stabilizing controllers from input-output data in the continuous-time, linear, and time-invariant domain. 
Rather than relying on measurements or reliable estimates of input and output time derivatives, the proposed approach uses filters to derive a parameterization of the system dynamics. 
This parameterization is amenable to the application of linear matrix inequalities enabling the design of stabilizing output feedback controllers from input-output data and the knowledge of the order of the system.
\end{abstract}


\section{Introduction\label{sec:intro}}

This letter aims to extend the results given in \cite{11045686} to the case where only input-output measurements are available.
To pursue this objective, some filters, borrowed from the adaptive control literature, are first used to obtain a state-space system representation
of the dynamics of the plant so that stabilization of such a system ensures stabilization of the closed loop. Although the obtained system is linear and time-invariant (LTI), it is affected by a disturbance whose
dynamics depend just on those of the filters. Hence, the framework proposed in \cite{11045686} is adapted to handle this disturbance.

As for the framework proposed in \cite{11045686}, and unlike most existing techniques  \cite{berberich2021data,lee2022formulations,bisoffi2022data,bisoffi2022learning,10384008,rapisarda2023fundamental,chu2023data,10623295,10323524,lopez2024data,10268593,10886882,schmitz2024continuous},
the approach proposed in this letter does not assume that measurements or reliable estimates of the successive time derivatives of the output are available.
Other approaches to designing stabilizing controllers from input-output data
that do not use the time derivatives of input and output have been proposed in \cite{2410.24167,2505.22505} for the single-input and single-output (SISO) and for the multi-input and multi-output 
(MIMO) cases, respectively.
The main differences from \cite{2410.24167,2505.22505} are the following:
\begin{enumerate*}[label*=(\roman*)]
	\item the data-driven parameterization proposed in this letter is independent of the time derivatives of the filter state;
	\item in the SISO case, the data matrix in this letter has $2n+1$ rows rather than $3n+1$ as in \cite{2410.24167,2505.22505};
	\item differently from  \cite{2410.24167,2505.22505}, where a dynamic extension is used to deal with unmeasurable states, in this letter a filtering strategy is proposed to remove disturbances;
	\item  sufficient conditions ensuring that the data matrix has full rank are provided.
\end{enumerate*}

\subsubsection*{Notation}
Let $I_n$ and $0_{n,m}$ denote the $n$-dimensional \emph{identity} and $(n\times m)$-dimensional \emph{zero} matrices, respectively. 
Let $\bbe{i}{n}$ denote the $i$th column of $I_n$.
Let \[\swap_n  =[\begin{array}{ccc}
	\bbe{n}{n} & \cdots & \bbe{1}{n}
\end{array}]\]
denote the $n$-dimensional \emph{reversal} matrix.
Let $\col(A,B)=[\begin{array}{cc}
	A^\top & B^\top
\end{array}]^\top$.
Let $\lambda(A)$ denote the \emph{spectrum} of $A$.
The symbol $\otimes$ denotes the Kronecker product.
A matrix is \emph{Toeplitz} if each descending diagonal from left to right is constant. 
A polynomial is \emph{Hurwitz} if all its roots have negative real parts.
A matrix $A\in\R^{n\times n}$ is \emph{Hurwitz} if its characteristic polynomial is Hurwitz.
Two polynomials in a single variable are \emph{coprime} if they do not have a common root.
Let $A^\dagger$ be the \emph{Moore-Penrose pseudoinverse} of the matrix $A$.
$\R[c_1,\dots,c_n]$ denotes the \emph{ring} of all the polynomials in $c_1,\dots,c_n$ with coefficients in $\R$.
Let $\imath=\sqrt{-1}$ denote the \emph{imaginary unit} in $\C$.
Let $\He(A)=A+A^\top$ and $\Sk(A)=A-A^\top$ be the \emph{symmetric} and \emph{skew-symmetric} parts of $A$, respectively.
The matrix $A\in\R^{n\times n}$ is \emph{positive definite}, denoted $A\succ 0$, if $x^\top A x>0$, $\forall x\in\R^n\setminus\{0\}$.
Let $\han{n}(\{d_k\}_{k=1}^{N-1})$ denote the \emph{Hankel matrix of depth $n$} of the discrete-time sequence $\{d_k\}_{k=1}^{N-1}$.
The sequence $\{d_k\}_{k=1}^{N-1}$, $d_k\in\R^m$, is \emph{persistently exciting of order $n$} if $\rank(\han{n}(\{d_k\}_{k=1}^{N-1}))=nm$.

\section{The proposed data-driven approach for the synthesis of an output feedback controller}\label{sec:pseudocode}
Consider a continuous-time, LTI, MIMO system described in differential operator representation \cite[Sec.~2.1]{kailath1980linear}
\begin{multline}\label{eq:diffEqMIMO}
	y^{(n)}(t)+A_{1}y^{(n-1)}(t)+\cdots+A_{n-1}y^{(1)}(t)+A_ny(t)\\
	=B_{1}u^{(n-1)}(t)+\cdots+B_{n-1}u^{(1)}(t)+B_nu(t),
\end{multline}
where $y(t)\in\R^p$, $u(t)\in\R^m$, $A_i\in\R^{p\times p}$, and $B_i\in\R^{p\times m}$, $i=1,\dots,n$.
Assume that a single input–output trajectory of system~\eqref{eq:diffEqMIMO}
in the interval $[0,T]$ is available:
\begin{equation*}
	\mD =\{(u(t),y(t))\in\R^m\times \R^p, t\in[0,T],\text{such that }\eqref{eq:diffEqMIMO}\text{ holds}\}.
\end{equation*}
The main goal of this letter is to design a stabilizing output feedback controller for system~\eqref{eq:diffEqMIMO}
just using the dataset $\mD$ and the knowledge of the order $n$ of system~\eqref{eq:diffEqMIMO}. 
Note that the knowledge of $n$ is usually assumed when applying adaptive approaches for
the synthesis of an output feedback controller; see, e.g., \cite[Chap.~7]{ioannou2012robust}.
\Cref{alg:proceure} illustrates the proposed approach to
design a stabilizing output feedback controller by using the dataset $\mD$ and the knowledge of $n$.

\begin{algorithm}[htb!]
	\caption{The proposed data-driven approach}\label{alg:proceure}
	\begin{algorithmic}[1]
		\REQUIRE dataset $\mD$, parameters $c_1,\dots,c_n,\beta\in\R$,  order $n$ of~\eqref{eq:diffEqMIMO} such that
			any of its minimal realizations has dimension $p n$
		\ENSURE a stabilizing output feedback controller or a \textbf{failure}
		\STATE define the matrices 
		\begin{gather}\label{eq:REF}
			\AREF = \left[
			\begin{array}{c}
				\begin{array}{c|c}
					0_{n-1,1} & I_{n-1}
				\end{array}\\
				\hline
				\begin{array}{ccccc}
					-c_n & \cdots & -c_1 
				\end{array}
			\end{array}\right],  \;\BREF = \left[\begin{array}{c}
				0_{n-1,1}\\
				\hline
				1
			\end{array}\right],
		\end{gather}
		\vspace{-2ex}
		\STATE let $\Au=\AREF\otimes I_m$, $\Bu=\BREF\otimes I_m$, $\Ay=\AREF\otimes I_p$, and $\By=\AREF\otimes I_p$
		\STATE filter the data in $\mD$ using the following filters
		\begin{subequations}
			\begin{align}
				\dot{\zeta}(t) & = \Au \zeta(t) + \Bu u(t), & \zeta(0)&=0_{mn,1},\label{eq:filtu}\\
				\dot{\mu}(t) & = \Ay \mu(t) + \By y(t), & \mu(0)&=0_{pn,1},\label{eq:filty}\\
				\dot{\phi}(t) & = -\beta \phi(t)+\chi(t), & \phi(0) & = 0_{2n,1},\label{eq:F1}\\
				\dot{\upsilon}(t) & = -\beta \upsilon(t) + u(t), & \upsilon(0)&=0,\label{eq:F2}\\
				\delta(t)  &= \chi(t)-\beta \phi(t),\label{eq:F3}
			\end{align}
			\label{eq:MIMOfilters}%
		\end{subequations}
		where $\chi(t)=\col(
			\zeta(t) , \mu(t)
		)$
		\STATE fix sampling times 
		$0<t_1<t_2\dots<t_N$  and define 
		\begin{subequations}
			\begin{align}
				\Delta_N & =  \left[\begin{array}{ccc}
					\delta(t_1) & \dots & \delta (t_N)
				\end{array}\right]\in\R^{(m+p)n\times N},\\
				\Phi_N & =  \left[\begin{array}{ccc}
					\phi (t_1) & \dots & \phi (t_N)
				\end{array}\right]\in\R^{(m+p)n\times N},\\
				\Upsilon_N & =  \left[\begin{array}{ccc}
					\upsilon (t_1) & \dots & \upsilon (t_N)
				\end{array}\right]\in\R^{m\times N}
			\end{align}
			\label{eq:datamatr}%
		\end{subequations}
		\vspace{-3ex}
		\STATE find a matrix $W_N$ such that,  for all $x_0\in\R^{pn}$, $[\begin{smallmatrix}
			 (\mathe^{\Ay^\top t_1}-\mathe^{-\beta t_1} I_{pn})x_0 & \cdots &
			  (\mathe^{\Ay^\top t_N}-\mathe^{-\beta t_N} I_{pn})x_0
		\end{smallmatrix}]W_N=0$
		\STATE define the filtered data
		\begin{align}\label{eq:filtered}
			\bar{\Delta}_N&=\Delta_N W_N,&
			\bar{\Phi}_N&=\Phi_NW_N,&
			\bar{\Upsilon}_N&=\Upsilon_NW_N
		\end{align}
		\vspace{-3ex}
		\IF{$\rank(\col(
				\bar{\Phi}_N,
				\bar{\Upsilon}_N)) = (m+p)n + m$}
			\STATE solve the LMI
			\vspace{-1ex}
			\begin{subequations}
					\begin{gather}
					\He(\bar{\Delta}_N  Z^\top)\prec 0,\qquad 
					\He(Z \bar{\Phi}_N^\top )  \succ 0,\\
					\Sk(Z \bar{\Phi}_N^\top )  = 0_{2n,2n}
				\end{gather}
				\label{eq:LMI}%
			\end{subequations}
			\vspace{-3ex}
			\STATE compute the gain 
				\begin{equation}\label{eq:Kstab}
					\KF  = \bar{\Upsilon}_N   Z^\top(Z\bar{ \Phi}_N^\top )^{-1}
				\end{equation}
				\vspace{-3ex}
			\RETURN the controller obtained by interconnecting~\eqref{eq:filtu},~\eqref{eq:filty}
				and letting $u(t)=\KF\col(
					\zeta(t) , \mu(t))$
		\ELSE
			\STATE \textbf{failure} the proposed approach cannot be applied
		\ENDIF
	\end{algorithmic}
\end{algorithm}

The effectiveness of \Cref{alg:proceure}  
is proven in \Cref{sec:SISO} in the SISO case and in \Cref{sec:MIMO} for a class of MIMO systems.

\section{The single-input single-output case}\label{sec:SISO}
The main goal of this section is to prove the effectiveness of \Cref{alg:proceure} to solve the problem of designing
a stabilizing output feedback controller in the SISO case. First, in \Cref{sec:inputoutput}, it is shown how to obtain
a model whose state is measurable from input-output data. Then, in \Cref{sec:synthesis}, it is shown how the gathered measurements
can be used to design a stabilizing output feedback controller.

\subsection{Data-driven modeling with input--output data}\label{sec:inputoutput}

If $p=m=1$ (i.e., in the SISO case), then~\eqref{eq:diffEqMIMO} becomes
\begin{multline}\label{eq:diffEq}
	y^{(n)}(t)+a_{1}y^{(n-1)}(t)+\cdots+a_{n-1}y^{(1)}(t)+a_ny(t)\\
	=b_{1}u^{(n-1)}(t)+\cdots+b_{n-1}u^{(1)}(t)+b_nu(t),
\end{multline}
where $y(t)\in\R$, $u(t)\in\R$, $y^{(i)}(t)=\frac{\D^i y(t)}{\D t^i}$, 
$u^{(i)}(t)=\frac{\D^i u(t)}{\D t^i}$, $a_i\in\R$, and $b_i\in\R$, $i=0,\dots,n$.
The following \namecref{ass:coprime} is made hereafter in this section.

\begin{asss}
	\label{ass:coprime}
	The polynomials $s^n+a_1s^{n-1}+\dots+a_n$ and $b_1s^{n-1}+b_2s^{n-2}+\dots+b_n$ are coprime.
\end{asss}

	\Cref{ass:coprime} entails the fact that the transfer function \[\frac{b_1s^{n-1}+b_2s^{n-2}+\dots+b_n}{s^n+a_1s^{n-1}+\dots+a_n}\]
	completely describes system~\eqref{eq:diffEq}. 

The description~\eqref{eq:diffEq}  is equivalent to the state-space form
	\begin{align}
		\label{eq:systemIO}%
		\dot{x}(t) & = Ax(t)+Bu(t),&
		y(t) & = Cx(t),
	\end{align}
where $x(t)\in\R^n$,
$A\in\R^{n\times n}$, $B\in\R^{n}$, and $C\in\R^{1\times n}$ are in \emph{observability canonical form}, i.e.,
$C=[\begin{array}{cc}
	0_{1,n-1} & 1
\end{array}]$, 
\begin{align*}
A  &= \left[\begin{array}{c|c}
	\begin{array}{c}
		0_{1,n-1}\\
		\hline
		I_{n-1}
	\end{array} 
	& \begin{array}{c}
		-a_n\\
		\vdots \\
		-a_1
	\end{array}
\end{array}\right], & 
B&=\left[\begin{array}{c}
	b_{n}\\
	\vdots \\ 
	b_1
\end{array}\right],
\end{align*}
by  \cite[Sec.~2.3]{kailath1980linear}. In particular, by letting
\begin{equation*}
	O  = \left[\begin{array}{c}
		C\\
		\vdots \\
		CA^{n-1}
	\end{array}\right]\!\!, \;
	M  = \left[\begin{array}{ccccccc}
		0 & \cdots & 0& 0\\
		\vdots & \vdots & \ddots &\vdots  \\
		CA^{n-2}B & \cdots & CB & 0
	\end{array}\right]\!\!,
\end{equation*}
be the observability matrix and the Toeplitz matrix of the Markov parameters
of system~\eqref{eq:systemIO}, the relation between the initial conditions
of the descriptions~\eqref{eq:diffEq} and~\eqref{eq:systemIO} is
\begin{equation}\label{eq:initialrelation}
	\left[\begin{array}{c}
		y^{(0)}(0)\\
		\vdots\\
		y^{(n-1)}(0)
	\end{array}\right]=O x(0) + M \left[\begin{array}{c}
		u^{(0)}(0)\\
		\vdots\\
		u^{(n-1)}(0)
	\end{array}\right].
\end{equation}

By \cite[Sec.~2.3 and Sec.~5.1]{kailath1980linear}, under \Cref{ass:coprime}, the pair $(A,B)$ is reachable and the pair $(A,C)$ is observable.

Since the time derivatives of the output $y(\cdot)$ and of the input $u(\cdot)$ are usually not measurable in practical settings, 
consider the filters~\eqref{eq:filtu} and~\eqref{eq:filty} of the input and output of~\eqref{eq:systemIO}.
\begin{remma}\label{rem:adaptive}
	The filters~\eqref{eq:filtu} and~\eqref{eq:filty} are usually employed in the adaptive control literature \cite{sastry2011adaptive,ioannou2012robust} to estimate the parameter vector $\theta=[\begin{array}{cccccc}
		b_1 & \cdots & b_n & a_1 & \cdots & a_n
	\end{array}]^\top$. In fact,  
	letting $\chi(t)=\col(
		\zeta(t) , \mu(t)
	)$ and $\psi(t)=y(t)+[\begin{array}{ccccccc}
	-c_n & -c_{n-1} & \cdots & -c_1
	\end{array}]\mu(t)$, the vector $\theta$ can be estimated using the \emph{gradient algorithm}
	\begin{equation*}
		\dot{\hat{\theta}}(t) = \frac{\psi(t)-\hat{\theta}^\top(t)\chi(t)}{1+\chi^\top(t)\chi(t)}\chi(t).
	\end{equation*}
	If $\chi(t)$ is \emph{persistently exciting} according to \cite[Def.~4.3.1]{ioannou2012robust}, then $\hat{\theta}(t)$ converges exponentially to $\theta$; see \cite[Thm.~4.3.2]{ioannou2012robust}.
\end{remma}

Unlike adaptive approaches, which estimate the parameters of system~\eqref{eq:diffEq}
 to design a control law, the main goal of this letter is to design a stabilizing control law using only input-output data,
bypassing the identification step. To this end, consider the next \namecref{lem:interconnection}, proved in Appendix~\ref{app:proof2}.

\begin{theom}
	\label{lem:interconnection}
	Let \Cref{ass:coprime} hold and  consider the interconnection of systems~\eqref{eq:systemIO},~\eqref{eq:filtu}, and~\eqref{eq:filty}, whose dynamics are
	\begin{align}
		\label{eq:intercon}%
		\dot{\eta}(t) &= \AIN \eta(t)+\BIN u(t), &
		\chi(t) &=\CIN \eta(t),
	\end{align}
	where $\eta(t) =\col(
		\zeta(t) , \mu(t) , x(t)
	)$, $\eta(t)\in\R^{3n}$, $\chi(t)=\col(
		\zeta(t) , \mu(t) 
	)$, $\chi(t)\in\R^{2n}$, $\CIN =\left[\begin{array}{ccc}
		I_{2n} & 0_{2n,n}
	\end{array}\right]$, and
	\begin{subequations}
		\begin{align}
			\AIN & = \left[\begin{array}{ccc}
				\AREF & 0_{n,n} & 0_{n,n}\\
				0_{n,n} & \AREF & \BREF C\\
				0_{n,n} & 0_{n,n} & A
			\end{array}\right], & \BIN & = \left[\begin{array}{c}
				\BREF \\ 0_{n,1} \\ B
			\end{array}\right].
		\end{align}
		\label{eq:intermatrix}%
	\end{subequations}
	By letting $x(0)=x_0$,
	the dynamics of $\chi(t)$ are given by
	\begin{equation}\label{eq:dynIO}
		\dot{\chi}(t)=\AF \chi(t) + \BF u(t)+ \GF \mathe^{\AREF^\top t}x_0,
	\end{equation}
	the pair $(\AF,\BF)$ is reachable, and
	\begin{gather*}
		\AF  = \left[\begin{array}{c|c}
			\begin{array}{c}
			\begin{array}{c|c}
				0_{n-1,1} & I_{n-1}
			\end{array}\\
			\hline
			\begin{array}{ccc}
				-c_n & \cdots & -c_1
			\end{array}
		\end{array}
		& 0_{n,n}\\
		\hline
		\begin{array}{c}
			0_{n-1,n}\\
			\hline
			\begin{array}{ccc}
				b_n & \cdots & b_1
			\end{array}
		\end{array}
		&
		\begin{array}{c}
			\begin{array}{c|c}
				0_{n-1,1} & I_{n-1}
			\end{array}\\
			\hline
			\begin{array}{ccc}
				-a_n & \cdots & -a_1
			\end{array}
		\end{array}
		\end{array}\right],\\
		\BF = \left[\begin{array}{c}
			0_{n-1,1}\\ \hline 1 \\ \hline 0_{n,1}
		\end{array}\right],\qquad
		\GF = \left[
		\begin{array}{c}
			0_{2n-1,n}\\
			\hline
			\begin{array}{c|c}
				0_{1,n-1} & 1
			\end{array}
		\end{array}\right].
	\end{gather*}
\end{theom}

By \Cref{lem:interconnection}, the dynamics of the state  $\chi(t)\in\R^{2n}$ of the filters~\eqref{eq:filtu} and~\eqref{eq:filty} are LTI and are affected
by a disturbance that depends on $c_1,\dots,c_n$ and on the initial condition $x_0$
of system~\eqref{eq:systemIO}. Consider the next \namecref{prp:stabilizing}.

\begin{crl}
	\label{prp:stabilizing}
	Let \Cref{ass:coprime} hold and let $\KF$ be such that $\AF+\BF\KF$ is Hurwitz.
	If the polynomial $s^n+c_1 s^{n-1}+\cdots+c_{n-1}s+c_n$ is Hurwitz, then the origin is globally asymptotically stable for the feedback interconnections
	of system~\eqref{eq:systemIO}, the filters~\eqref{eq:filtu} and~\eqref{eq:filty}, and $u(t)=\KF \chi(t)$.
\end{crl}

\begin{proof}
	The proof follows directly from Appendix~\ref{app:proof2} and
	the fact that if  $s^n+c_1 s^{n-1}+\cdots+c_{n-1}s+c_n$ is Hurwitz and $\KF$ is such that $\AF+\BF\KF$ is Hurwitz, then the matrix 
	\[\left[\begin{array}{cc}
		\AF +\BF \KF & \GF \\
		0_{n,2n} & \AREF^\top
	\end{array}\right]\] is Hurwitz.
\end{proof}

By \Cref{prp:stabilizing}, if $s^n+c_1 s^{n-1}+\cdots+c_{n-1}s+c_n$ is Hurwitz, then it suffices to
design a feedback gain $\KF$ such that $\AF+\BF\KF$ is Hurwitz 
yielding a stabilizing controller consisting of the filters~\eqref{eq:filtu},~\eqref{eq:filty} and the control law $u(t)=\KF \chi(t)$. 
Since the dynamics of $\chi$ are LTI with a disturbance generated by
an LTI known exosystem, this objective can be pursued by suitably adapting the framework proposed in \cite{11045686}.
Namely, consider the additional filters~\eqref{eq:F1},~\eqref{eq:F2}, and~\eqref{eq:F3},
where $\beta\in\R$ is a parameter assumed to satisfy the following \namecref{ass:noresonance}.
 
\begin{asss}
	\label{ass:noresonance}
	The parameter $\beta$ satisfies $-\beta\notin\lambda(\AREF)$.
\end{asss}

\Cref{ass:noresonance} holds for all $\beta\in\R$ such that $(-\beta)^n+c_1(-\beta)^{n-1}+\cdots+c_n\neq 0$.
Consider the following \namecref{lem:relation}, whose proof is given in Appendix~\ref{app:proof4}.

\begin{lemm}
	\label{lem:relation}
	Let \Cref{ass:coprime,ass:noresonance} hold. By letting $\Gamma^\star\in\R^{2n\times n}$ be the unique solution to 
	$\Gamma \AREF^\top + \beta \Gamma  = \GF$,
	define $\epsilon(t)=\Gamma^\star (\mathe^{\AREF^\top t}-\mathe^{-\beta t} I_{n})x_0$.
	For all $t\geq 0$, one has
	\begin{equation}\label{eq:deltaAF}
		\delta(t) = \AF \phi(t) + \BF \upsilon(t) + \epsilon(t).
	\end{equation}
\end{lemm}

Building on \Cref{lem:relation} and using a construction similar to that employed in \cite{8933093,11045686}, 
 consider the following \namecref{ass:fullrank}.

\begin{asss}\label{ass:fullrank}
    Let $E_N=[\begin{array}{ccc}
    	\epsilon(t_1) & \cdots & \epsilon(t_N)
    \end{array}]\in\R^{2n\times N}$.
	There is a matrix $W_N\in\R^{N\times \bar{N}}$, $\bar{N}\in\N$, such that $E_N W_N=0_{2n,\bar{N}}$, and, by defining the filtered data matrices as in~\eqref{eq:filtered},
	the following rank condition holds
	\begin{equation}\label{eq:rank}
		\rank(\col(\bar{\Phi}_N,\bar{\Upsilon}_N))
		= 2n + 1.
	\end{equation}
\end{asss}
By \Cref{lem:relation}, if \Cref{ass:coprime,ass:noresonance,ass:fullrank} hold, then the matrices $\AF$ and $\BF$ can be obtained as
\begin{equation}\label{eq:estimation}
	[\begin{array}{cc}
		\AF & \BF
	\end{array}] = \bar{\Delta}_N (\col(\bar{\Phi}_N,\bar{\Upsilon}_N))^\dagger.
\end{equation}
In fact,  by \Cref{lem:relation}, one has $\Delta_N=\AF \Phi_N+\BF\Upsilon_N+E_N$.
Hence, since $W_N$ is such that $E_NW_N=0$, one has $\bar{\Delta}_N=\AF \bar{\Phi}_N+\BF\bar{\Upsilon}_N$.
Hence, in principle, a gain $\KF$ such that $\AF+\BF\KF$ is Hurwitz can be designed by estimating $\AF$ and $\BF$ using~\eqref{eq:estimation},
provided that \Cref{ass:coprime,ass:noresonance,ass:fullrank} hold.
In \Cref{sec:synthesis}, it is shown how to directly design $\KF$ using the matrices in~\eqref{eq:filtered}
bypassing the intermediate estimation step. 

The main goal of the remainder of this section is to provide conditions ensuring that \Cref{ass:fullrank} is satisfied.
Consider the following \namecref{lem:null}, whose proof is given in Appendix~\ref{app:proof5}.

\begin{lemm}
	\label{lem:null}
	Let $N\geq n+1$ and $t_1<\dots<t_{N}$ be fixed.
	There exist $N-n$ linearly independent $\alpha\in\R^{N}$ such that
	\[\sum_{i=1}^{N}\alpha_{i} (\mathe^{\AREF^\top t_i}-\mathe^{-\beta t_i} I_{n})=0_{n,n}.\]
\end{lemm}

By \Cref{lem:null}, if $N\geq n+1$, then there exists $W_N\in\R^{N\times (N-n)}$, $\rank(W_N)=N-n$,
depending only on $t_1\dots,\dots,t_{N}$, $c_1,\dots,c_n$, and $\beta$, such that
$E_N W_N=0_{2n,N-n}$, $\forall x_0\in\R^n$. 
In particular, the matrix $W_N$ can be obtained by stacking the $N-n$ linearly independent $\alpha\in\R^{N}$ characterized in
\Cref{lem:null}. 
The next \namecref{ex:uniform} illustrates how to obtain $\bar{\Delta}_N$, $\bar{\Phi}_N$, and $\bar{\Upsilon}_N $
using a discrete-time filter in the case that
$t_1<\dots < t_N$ are uniformly spaced.

\begin{remma}\label{ex:uniform}
	Let $N\geq n+2$.
	If $t_1,\dots,t_N$ are uniformly spaced, $t_i = i \TS$, $i=1,\dots,N$, for some \emph{sampling time} $T_S>0$, then
	by the Cayley-Hamilton theorem \cite[Ex.~7.2.2]{meyer2023matrix}, by letting $W_N \in\R^{N\times (N-n-1)}$ 
	be the lower triangular Toeplitz matrix whose first column is $[\begin{array}{cccccc}
		w_{n+1} &  \cdots & w_{1} & w_{0} & 0_{1,N-n-2}
	\end{array}]^\top$,  where $w_0z^{n+1}+\cdots+w_{n+1}=(z-\mathe^{-\beta \TS})\det(zI-\mathe^{\AREF^\top \TS})$,
	one has  $\rank(W_N)=N-n-1$ and $E_NW_N=0$, $\forall x_0\in\R^n$.
	Note that $\bar{\Delta}_N $, $\bar{\Phi}_N $, and $\bar{\Upsilon}_N $ can be obtained by feeding the discrete-time signals $\{\delta_i(t_k)\}_{k=1}^N$,
	$\{\phi_i(t_k)\}_{k=1}^N$, $i=1,\dots,2n$, and $\{\upsilon(t_k)\}_{k=1}^N$, to the finite impulse response (FIR) filter 
	\[\frac{z^{n+1}+w_1z^n+\cdots+w_nz+w_{n+1}}{z^{n+1}},\]
	and discarding its first $n$ outputs.
\end{remma}

Consider the following \namecref{ass:noresonance2}.

\begin{asss}
	\label{ass:noresonance2}
	The polynomials $s^n+a_1s^{n-1}+\dots+a_n$ and $s^n+c_1s^{n-1}+\dots+c_n$  are coprime.
\end{asss}

By letting $a_1,\dots,a_n\in\R$ and $\beta\in\R$ be fixed,
 let $q\in\R[c_1,\dots,c_n]$ be the resultant \cite[\S6, Chap.~3]{cox2013ideals} of
$(s^n+a_1s^{n-1}+\dots+a_n)(s+\beta)$ and $s^n+c_1s^{n-1}+\dots+c_n$.
\Cref{ass:noresonance,ass:noresonance2} hold if and only if $[\begin{array}{ccc}
		c_1 & \cdots & c_n
	\end{array}]^\top\notin\mathbf{V}(q)$, where $\mathbf{V}(q)=\{[\begin{array}{ccc}
	c_1 & \cdots & c_n
	\end{array}]^\top\in\R^n:q(c_1,\dots,c_n)=0\}$.
Thus, \Cref{ass:noresonance,ass:noresonance2} hold generically.
These assumptions are made in \Cref{thm:pe}, proved in Appendix~\ref{app:proof6}, to simplify the analysis of the interconnection of~\eqref{eq:systemIO} and \eqref{eq:MIMOfilters},
ruling out the presence of resonances, to establish sufficient conditions ensuring that \Cref{ass:fullrank} holds.
Since the interconnection of~\eqref{eq:diffEq} and~\eqref{eq:MIMOfilters} is not reachable (see Appendix~\ref{app:proof6}), such a \namecref{thm:pe}
does not follow directly from the Willems et al. fundamental lemma \cite[Cor.~2]{willems2005note} or \cite[Lem.~1]{lopez2022continuous}.

\begin{prpp}
	\label{thm:pe}
	Let \Cref{ass:coprime,ass:noresonance,ass:noresonance2} hold.
	Let $\TS>0$ be a given sampling time such that 
	\begin{multline}
		\label{eq:nonpato}
		\nexists \ell\in\{-\beta\}\cup \lambda(A)\cup\lambda(\AREF), h\in\mathbb{Z}\setminus\{0\} \text{ such that }\\\ell+\imath {2h\pi}\TS^{-1}\in\{-\beta\}\cup \lambda(A)\cup\lambda(\AREF),
	\end{multline}
	By letting $N\geq 8n+4$, $t_k=k\TS$, $k\in\{0,\dots,N\}$, let $W_N\in\R^{N\times (N-n-1)}$ 
	be the lower triangular Toeplitz matrix whose first column is $[\begin{array}{cccccc}
		w_{n+1} &  \cdots & w_{1} & w_{0} & 0_{1,N-n-2}
	\end{array}]^\top$,  where $w_0z^{n+1}+w_1z^n+\cdots+w_{n+1}=(z-\mathe^{-\beta \TS})\det(zI-\mathe^{\AREF^\top \TS})$. 
	By applying to the interconnection of system~\eqref{eq:systemIO} 
	and the filters~\eqref{eq:MIMOfilters},	the input 
	 $u(t)=d_k$, $\forall t\in[k\TS,(k+1)\TS)$, $k\in\{0,\dots,N-1\}$, if $\{d_k\}_{k=1}^{N-1}$ is persistently exciting of order $4n+2$, then \Cref{ass:fullrank} holds. 
\end{prpp}

Note that~\eqref{eq:nonpato} holds for almost all $\TS>0$ (see the discussion in \cite{11045686}). 
Therefore, since \Cref{ass:noresonance,ass:noresonance2} are also generically satisfied, \Cref{thm:pe} provides an easily implementable approach to gather informative data.

By letting  $\LF\in\R^{1\times 2n}$, the result of \Cref{thm:pe} can be extended to include a feedback term in the control input:
	$u(t) = \LF \chi(t) + d_k$, $\forall t\in[k\TS,(k+1)\TS)$.
If \Cref{ass:coprime,ass:noresonance} hold, $\TS>0$ is such that~\eqref{eq:nonpato} holds with $\lambda(A)\cup\lambda(\AREF)$ substituted by $\lambda(\AF+\BF\LF)$,
$\lambda(A+\BF\LF)\cap \lambda(\AREF)=\emptyset$, and $\{d_k\}_{k=1}^{N-1}$ is persistently exciting of order $4n+2$, then a reasoning wholly similar to that given in Appendix~\ref{app:proof6}
can be used to conclude that \Cref{ass:fullrank} holds.

\subsection{Data-driven synthesis of a stabilizing gain\label{sec:synthesis}}

By \Cref{prp:stabilizing}, if $\KF$ is such that $\AF+\BF\KF$ is Hurwitz, then the feedback interconnection of system~\eqref{eq:systemIO}, the filters~\eqref{eq:filtu},~\eqref{eq:filty},
and $u(t)=\KF \chi(t)$ is globally asymptotically stable. Thus, consider the following \namecref{thm:stabil}, which follows directly from~\eqref{eq:estimation} and \cite[Thm.~2]{11045686}.

\begin{theom}\label{thm:stabil}
	Let \Cref{ass:coprime,ass:noresonance,ass:fullrank} hold. Then, any matrix $Z\in\R^{2n\times h}$ such that the LMI~\eqref{eq:LMI} holds
	is such that the gain $\KF$ given in~\eqref{eq:Kstab}
	is stabilizing. Conversely, if $\KF$ is stabilizing, then it can be written as in~\eqref{eq:Kstab}, with $Z$ satisfying~\eqref{eq:LMI}.
\end{theom}

It is worth pointing out that an approach wholly similar to that given in \cite[Thm.~3]{11045686}
can be used to guarantee robustness of such a scheme with respect to  measurement and input noises.

\section{Extension to a class of multi-input multi-output systems}\label{sec:MIMO}
The main objective  of this section is to extend the results given in \Cref{sec:SISO} to a class of MIMO systems.
Toward this goal, consider a generic MIMO system in state space form
\begin{align}
	\label{eq:sysMIMO}
	\dot{\bar{x}}(t) & = \bar{A} \bar{x}(t)+\bar{B} u(t), & y(t)=\bar{C} \bar{x}(t),
\end{align}
with $\bar{x}(t)\in\R^h$, $u(t)\in\R^m$, and $y(t)\in\R^p$.  The next \namecref{ass:MIMOass} is made throughout this section to ensure that system~\eqref{eq:sysMIMO}
is a realization of the differential form~\eqref{eq:diffEqMIMO}.
\begin{taggedassumption}{\ref{ass:coprime}'}
	\label{ass:MIMOass}
	System~\eqref{eq:sysMIMO} is reachable and observable, $p n = h$, and
	$\rank(\bar{O})=h$, $\bar{O}=\col(\bar{C},\bar{C}\bar{A},\cdots,\bar{C}\bar{A}^{n-1})$.
\end{taggedassumption}
\begin{remma}\label{rem:comments}
	Under \Cref{ass:MIMOass}, 
	by letting $\bar{M}$ be the matrix obtained from $M$ by substituting $A$, $B$, and $C$ with $\bar{A}$, $\bar{B}$, and $\bar{C}$, respectively, one has that \[y^{(n)}(t) = \bar{C} \bar{A}^n \bar{x}(t)+\sum_{i = 0}^{n-1}\bar{C}\bar{A}^{n-i-1}\bar{B} u^{(i)}(t)\] and $\bar{x}(t)$ can be uniquely expressed in terms of $y^{(0)}(t),\cdots,y^{(n-1)}(t),u^{(0)}(t),\cdots,u^{(n-1)}(t)$ as 
	\begin{multline*}
		\bar{x}(t)=O^{-1}(\col(y^{(0)}(t),\cdots,y^{(n-1)}(t))\\-\bar{M}\col(u^{(0)}(t),\cdots,u^{(n-1)}(t))).
	\end{multline*}
	Thus, under \Cref{ass:MIMOass}, Eq.~\eqref{eq:diffEqMIMO} completely describes the dynamics of system~\eqref{eq:sysMIMO}.
\end{remma}

Define the matrices 
$\CM=[\begin{array}{cc}
	0_{p,p(n-1)} & I_p
\end{array}]$ and
\[\AM  = \left[\begin{array}{c|c}
	\begin{array}{c}
		0_{p,p(n-1)}\\
		\hline
		I_{p(n-1)}
	\end{array} 
	& \begin{array}{c}
		-A_n\\
		\vdots \\
		-A_1
	\end{array}
\end{array}\right],\qquad 
\BM=\left[\begin{array}{c}
	B_{n}\\
	\vdots \\ 
	B_1
\end{array}\right].\]
Under \Cref{ass:MIMOass}, the pair $(\AM,\BM)$ is reachable since the triplet $(\AM,\BM,\CM)$ is a realization of~\eqref{eq:diffEqMIMO}, and all minimal  (i.e., reachable and observable) realizations have the same order \cite[Sec.~17.1]{hespanha2018linear}.
In particular, by letting 
\begin{align*}
	\OM  &= \left[\begin{array}{c}
		\CM\\
		\vdots \\
		\CM\AM^{n-1}
	\end{array}\right]\\\
	\MM & = \left[\begin{array}{ccccccc}
		0 & \cdots & 0& 0\\
		\vdots & \vdots & \ddots &\vdots  \\
		\CM\AM^{n-2}B & \cdots & \CM\BM & 0
	\end{array}\right]\!\!,
\end{align*}
the relation between the initial conditions
of the description in differential operator form~\eqref{eq:diffEqMIMO} and
the one in state space form
\begin{align*}
	\dot{x}(t) & = \AM x(t) + \BM u(t), &
	y(t) & = \CM x(t),
\end{align*}
is, by \Cref{rem:comments}, given by the following relation
\begin{equation}\label{eq:initialrelationMIMO}
	\left[\begin{array}{c}
		y^{(0)}(0)\\
		\vdots\\
		y^{(n-1)}(0)
	\end{array}\right]=\OM x(0) + \MM \left[\begin{array}{c}
		u^{(0)}(0)\\
		\vdots\\
		u^{(n-1)}(0)
	\end{array}\right].
\end{equation}
Hence, consider the following \namecref{thm:interMIMO}, 
whose proof is given in Appendix~\ref{app:proof7},
that is the equivalent of \Cref{lem:interconnection}
in the MIMO case considered in this section.

\begin{taggedtheorem}{\ref{lem:interconnection}'}\label{thm:interMIMO}
	Let \Cref{ass:MIMOass} hold and  consider the interconnection of systems~\eqref{eq:diffEqMIMO},~\eqref{eq:filtu},~\eqref{eq:filty}, whose dynamics are
	\begin{align}
		\label{eq:interconMIMO}%
		\dot{\eta}(t) &= \AMIN \eta(t)+\BMIN u(t), &
		\chi(t) &=\CMIN \eta(t),
	\end{align}
	where $\eta(t) =\col(
	\zeta(t) , \mu(t) , x(t)
	)$, $\eta(t)\in\R^{(m+2p)n}$, $\chi(t)=\col(
	\zeta(t) , \mu(t) 
	)$, $\chi(t)\in\R^{(m+p)n}$, $\CIN =\left[\begin{array}{ccc}
		I_{(m+p)n} & 0_{(m+p)n,pn}
	\end{array}\right]$, and
	\begin{subequations}
		\begin{align}
			\AMIN & = \left[\begin{array}{ccc}
				\Au & 0_{mn,pn} & 0_{mn,pn}\\
				0_{pn,mn} & \Ay & \By \CM\\
				0_{pn,mn} & 0_{pn,pn} & \AM
			\end{array}\right], \\ \BMIN & = \left[\begin{array}{c}
					\Bu \\ 0_{pn,m} \\ \BM
			\end{array}\right].
		\end{align}
		\label{eq:intermatrixMIMO}%
	\end{subequations}
	By letting $x(0)=x_0$, where $x(0)$ satisfies~\eqref{eq:initialrelationMIMO}, 
	the dynamics of $\chi(t)$ are given by
	\begin{equation}\label{eq:dynIOMIMO}
		\dot{\chi}(t)=\AEF \chi(t)+ \BEF u(t)+ \GEF \mathe^{\Ay^\top t}x_0,
	\end{equation}
	the pair $(\AEF,\BEF)$ is reachable, and
	\begin{gather*}
		\AEF  = \left[\begin{array}{c|c}
			\Au
			& 0_{mn,pn}\\
			\hline
			\begin{array}{c}
				0_{(n-1)p,nm}\\
				\hline
				\begin{array}{ccc}
					B_n & \cdots & B_1
				\end{array}
			\end{array}
			&
			\begin{array}{c}
				\begin{array}{c|c}
					0_{(n-1)p,p} & I_{(n-1)p}
				\end{array}\\
				\hline
				\begin{array}{ccc}
					-A_n & \cdots & -A_1
				\end{array}
			\end{array}
		\end{array}\right],\\
		\BEF = \left[\begin{array}{c}
			\Bu \\ 0_{pn,m}
		\end{array}\right],\quad \GEF = \left[\begin{array}{c}
		0_{mn,pn}\\ \By \CM
		\end{array}\right].
	\end{gather*}
\end{taggedtheorem}

By leveraging \Cref{thm:interMIMO}, the following \namecref{lem:relationMIMO}, which is the equivalent of \Cref{lem:relation}
in the MIMO case considered in this section and whose proof is given in Appendix~\ref{app:proof8}, characterizes the response of
the filters~\eqref{eq:F1},~\eqref{eq:F2}, and~\eqref{eq:F3}. 

\begin{taggedlemma}{\ref{lem:relation}'}
	\label{lem:relationMIMO}
	Let Assumptions~\ref{ass:MIMOass} and~\ref{ass:noresonance} hold. 
	Consider the response of the filters~\eqref{eq:F1},~\eqref{eq:F2}, and~\eqref{eq:F3}.  
	One has that
	\begin{equation}\label{eq:deltaAFrealization}
		\delta(t) = \AEF \phi(t) + \BEF \upsilon(t) + \epse(t),
	\end{equation}
	for all $t\geq 0$, where
	\[\epse(t)=\GEF(\Ay^\top+\beta I_{pn})^{-1}(\mathe^{\Ay^\top t}-\mathe^{-\beta t}I_{pn})x_0.\]
\end{taggedlemma}

Since the minimal polynomial of $\Ay$ is $s^n+c_1s^{n-1}+\dots+c_n$,
if the sampling times are uniformly spaced, i.e. $t_i=i\TS$, $i=1,\dots,N$, then the matrix
$W_N$ given in \Cref{ex:uniform} is such that
	$[\begin{array}{cccc}
		\epse(t_1) & \cdots & \epse(t_N)
	\end{array}]W_N=0$,
for all $x_0\in\R^{pn}$. Thus, consider the following \namecref{thm:peMIMO}, whose proof is given in Appendix~\ref{app:proof9},
that extends \Cref{thm:pe} to the MIMO case considered in this section.

\begin{taggedproposition}{\ref{thm:pe}'}\label{thm:peMIMO}
	Let Assumptions~\ref{ass:MIMOass} and~\ref{ass:noresonance} hold.
	Additionally, suppose that  $\sigma(\AM)\cap \sigma(\AREF)=\emptyset$.
	Let $\TS>0$ be a given sampling time such that 
	\begin{multline}
		\label{eq:nonpatoMIMO}
		\nexists \ell\in\{-\beta\}\cup \lambda(\AM)\cup\lambda(\AREF), h\in\mathbb{Z}\setminus\{0\} \text{ such that }\\\ell+\imath {2h\pi}\TS^{-1}\in\{-\beta\}\cup \lambda(\AM)\cup\lambda(\AREF),
	\end{multline}
	Let $W_N\in\R^{N\times (N-n-1)}$ 
	be the lower triangular Toeplitz matrix whose first column is $[\begin{array}{cccccc}
		w_{n+1} &  \cdots & w_{1} & w_{0} & 0_{1,N-n-2}
	\end{array}]^\top$,   $w_0z^{n+1}+\cdots+w_{n+1}=(z-\mathe^{-\beta \TS})\det(zI-\mathe^{\AREF^\top \TS})$. 
	By applying to the interconnection of system~\eqref{eq:diffEqMIMO} 
	and the filters~\eqref{eq:MIMOfilters},	the input 
	$u(t)=d_k$, $\forall t\in[k\TS,(k+1)\TS)$, $k\in\{0,\dots,N-1\}$, if $\{d_k\}_{k=1}^{N-1}$ is persistently exciting of order  $2((m+p)n+m)$, 
	then the matrices $\bar{\Delta}_N$, $\bar{\Phi}_N$, and $\bar{\Upsilon}_N$ defined in~\eqref{eq:filtered} satisfy
	\[\rank(\col(\bar{\Phi}_N,\bar{\Upsilon}_N))=(m+p)n+m.\]
\end{taggedproposition}

By the same reasoning given in \Cref{sec:synthesis}, under the hypotheses of \Cref{thm:peMIMO},
by solving the LMI~\eqref{eq:LMI}, a gain $\KEF$ such that $\AEF+\BEF\KEF$ is stabilizing can be obtained as in~\eqref{eq:Kstab}.
Then, by the same reasoning used to prove \Cref{prp:stabilizing}, 
a stabilizing controller can be obtained using the filters~\eqref{eq:filtu},~\eqref{eq:filty} and by letting $u(t)=\KEF \chi(t)$.

\section{Conclusions}\label{sec:conclusions}

A data-driven framework for the design of output feedback controllers for LTI continuous-time systems has been presented.
By combining input–output filtering with an LMI-based synthesis procedure, the proposed method avoids reliance on time derivatives and bypasses explicit model identification. 
Constructive conditions ensuring the applicability of the framework have been provided. 
Future work will focus on estimating the order $n$ of system~\eqref{eq:diffEqMIMO} from data and
on weakening \Cref{ass:MIMOass} in the MIMO case.

\appendix
\subsection{Proof of \Cref{lem:interconnection}\label{app:proof2}}
It is first proved that, $\forall k\in\N$, $k\geq 1$, one has
\begin{equation}
	\label{eq:AINkBIN}
	\AIN^k\BIN = \left[\begin{array}{c}
		\AREF^k\BREF \\
		\sum_{i=0}^{k-1}\AREF^i \BREF C A^{k-i-1}B\\
		A^k B\\
	\end{array}\right].
\end{equation}
By definition of $\AIN$ and $\BIN$, Eq.~\eqref{eq:AINkBIN} holds for $k=1$.
Assuming that it holds for some $k\in\N$, $k\geq 1$, by computing
$\AIN \cdot \AIN^k \BIN$, one has that~\eqref{eq:AINkBIN} holds with $k$ substituted by $k+1$.
Therefore, Eq.~\eqref{eq:AINkBIN} holds by induction.

It is then proved that 
\begin{subequations}
	\begin{gather}
				\AREF^k\BREF +\sum_{j=1}^k c_j \AREF^{k-j}\BREF=\bbe{n-k}{n}, k = 1,\dots,n-1,\label{eq:sum1}\\
		\AREF^{n+i}\BREF +\sum_{j=1}^{n}c_j \AREF^{n+i-j}\BREF=0_{n,1},  \forall i \in \N.\label{eq:sum2}
	\end{gather}
	\label{eq:sum}%
\end{subequations}
Note that $\BREF=\bbe{n}{n}$ and $\AREF \BREF + c_1 \BREF=\bbe{n-1}{n}$. Hence, Eq.~\eqref{eq:sum1} holds for $k=1$. Assume that~\eqref{eq:sum}
holds for some $k\in\N$, $1 \leq k < n$. Then, one has $\AREF \bbe{n-k}{n}=\bbe{n-k-1}{n}-c_{k+1}\bbe{n}{n}=
\AREF^{k+1}\BREF +c_1 \AREF^{k}\BREF+\cdots+c_k\AREF\BREF$. Hence, Eq.~\eqref{eq:sum1} holds by induction for
$k=1,\dots,n-1$. On the other hand, by the Cayley-Hamilton theorem, one has that~\eqref{eq:sum2} holds.

Equations~\eqref{eq:AINkBIN} and~\eqref{eq:sum} are now used to find a basis for the image of 
the reachability matrix $\RIN=[\begin{array}{ccc}
	\BIN & \cdots & \AIN^{3n-1}\BIN
\end{array}]$ of the pair $(\AIN,\BIN)$.
 By letting $c_0=1$, by~\eqref{eq:AINkBIN}, one has
 \begin{multline*}
 	\sum_{j=0}^k c_j \AIN^{k-j}\BIN \\
 	=\sum_{j=0}^k c_j \col\left(	\AREF^{k-j}\BREF , \sum_{i=0}^{k-j-1}\AREF^{k-j-i-1} \BREF C A^{i}B,A^{k-j} B\right).
 \end{multline*}
Note that 
\begin{multline*}
\sum_{j=0}^k c_j \sum_{i=0}^{k-j-1}\AREF^{k-j-i-1} \BREF C A^{i}B\\=\sum_{i=0}^{k-1}\left(\sum_{j=0}^{k-i-1}c_j \AREF^{k-j-i-1}\BREF\right) CA^iB.
\end{multline*}
By~\eqref{eq:sum1}, one has $\sum_{j=0}^{k-i-1}c_j \AREF^{k-j-i-1}\BREF= \bbe{n-k+i+1}{n}$, which implies, for $k=1,\dots,n-1$,
\begin{multline*}
\sum_{j=0}^k c_j \AIN^{k-j}\BIN \\
=\col\left(\bbe{n-k}{n},\allowbreak \sum_{i=0}^{k-1}\bbe{n-k+i+1}{n}CA^iB, \allowbreak\sum_{j=0}^k c_j A^{k-j} B\right).
\end{multline*}
Similarly, by~\eqref{eq:sum2}, one has, for $k=n,\dots,3n-1$,
\begin{multline*}
	\sum_{j=0}^n c_j \AIN^{k-j}\BIN \\= \col\left(0_{n,1},\allowbreak \sum_{i=0}^{n-1}\bbe{n-i}{n}CA^{k-i-1}B, \allowbreak \sum_{j=0}^n c_j A^{k-j} B\right).
\end{multline*}
Therefore, by letting $\HF\in\R^{3n\times 3n}$ be the upper triangular Toeplitz matrix 
whose first row is $[\begin{array}{cccccccc}
	1 & c_1 & \cdots & c_n & 0_{1,2n-1}
\end{array}]$, one has
\[\RIN \HF = \left[\begin{array}{cc}
		\swap_n & 0_{n,2n}\\
		M \swap_n & O R_{2n}\\
		F_n \swap_n & q_{\mathrm{r}}(A)R_{2n}
	\end{array}\right],\]
where $R_{2n} = [\begin{array}{cccc}
	B & AB & \cdots & A^{2n-1}B
\end{array}]$, $F_n=[\begin{array}{ccc}
	( \sum_{i = 0}^{n - 1} c_i A^{n  - i- 1} ) B & \cdots & B
\end{array}]$, and $q_{\mathrm{r}}(s)=s^n+c_1s^{n-1}+\cdots +c_n$.
Under \Cref{ass:coprime}, one has $\rank(R_{2n} )=n$.
Since $O$ is invertible and $\swap_n^{-1}=\swap_n$, by taking linear combinations of the columns of $\RIN \HF $, one has 
\begin{equation*}
	\SP(\RIN)=\SP\left(\left[\begin{array}{cc}
		I_n & 0_{n,n}\\
		0_{n,n} & I_{n}\\
		\mho_1 & \mho_2
	\end{array}\right]\right),
\end{equation*} 
where $\mho_1=F_n-q_{\mathrm{r}}(A)O^{-1}M$ and $\mho_2=q_{\mathrm{r}}(A)O^{-1}$.
 Since $(A,B,C)$ is a realization of~\eqref{eq:diffEq}, one has $CA^nO^{-1}=[\begin{array}{ccc}
 	-a_n & \cdots & -a_1
 \end{array}]$ and $[\begin{array}{ccc}
 	C A^{n-1}B & \cdots & CB
 \end{array}]-CA^{n}O^{-1}M=[\begin{array}{ccc}
 	b_n & \cdots & b_1
 \end{array}]$; see \Cref{rem:comments}. Hence, it results that
 \begin{align*}
 	C\mho_2 & =  \sum_{i=0}^{n-1}(c_{n-i}-a_{n-i}) CA^i O^{-1}\\
 	&=[\begin{array}{ccc}
 		c_n-a_n & \cdots & c_1-a_1
 	\end{array}].
 \end{align*}
Further, by \cite[Eq.~(28) at p.~324]{kailath1980linear}, one has 
\begin{align*}
C\mho_1&=CF_n-C\mho_2 M\\
&=[\begin{array}{cccc}
CA^{n - 1} B + \sum_{i = 1}^{n - 1} a_i CA^{n - 1 - i} B & \cdots 
& CB
\end{array}]\\
&=B^\top.
\end{align*} Finally, since $O^{-1}\BREF=\bbe{1}{n}$, by \cite[Eq.~(11) at p.~199]{kailath1980linear}, and $A^{i-1} \bbe{1}{n}=\bbe{i}{n}$, $i=1,\dots,n$,
 one has 
 \[\mho_2 \BREF=[\begin{array}{ccc}
 	c_n-a_n & \cdots & c_1-a_1
 \end{array}]^\top.\]
Thus, let 
\[P=\left[\begin{array}{ccc}
	I_n & 0_{n,n} & 0_{n,n}\\
	0_{n,n} & I_{n} & 0_{n,n}\\
	\mho_1 & \mho_2 & I_n
\end{array}\right].\]
By classical results on the Kalman decomposition for reachability \cite[Sec.~16.2]{hespanha2018linear}, one has
\begin{align*}
	P^{-1}\AIN P &= \left[\begin{array}{cc}
		\AINrr & \AINru\\
		0_{n,2n} & \AINuu
	\end{array}\right], & P^{-1}\BIN & = \left[\begin{array}{c}
		\BINr \\ 0_{n,1}
	\end{array}\right],
\end{align*}
and the pair $(\AINrr,\BINr)$ is reachable. Note that
$\CIN P=\CIN$, $\AREF+\BREF C \mho_2=A^\top$, $P^{-1}\col(\chi(t),x(t))=\col(\chi(t),x(t)-\mho_1\zeta(t)-\mho_2\mu(t) )$, 
$\AINrr  = \left[\begin{smallmatrix}
	\AREF & 0_{n,n}\\
	\BREF C \mho_1 & \AREF+\BREF C \mho_2
\end{smallmatrix}\right]$, $\BINr  = \left[\begin{smallmatrix}
\BREF \\ 0_{n,1}
\end{smallmatrix}\right]$, $\AINru = \left[\begin{smallmatrix}
0_{n,n}\\ \BREF C
\end{smallmatrix}\right]$, and $ \AINuu = A- \mho_2 \BREF C=\AREF^\top$.
Hence, the statement follows by the fact that $\AINrr  =\AF$, $\BINr=\BF$, $\AINru=\GF$,
and $\kappa(t)=x(t)-\mho_1\zeta(t)-\mho_2\mu(t)$ satisfies $\dot{\kappa}(t)=\AINuu \kappa(t)$
with $\kappa(0)=x_0$.
\subsection{Proof of \Cref{lem:relation}\label{app:proof4}}
Following a reasoning similar to that used in the proof of \cite[Lem.~1]{11045686}, the solution to system~\eqref{eq:F1},~\eqref{eq:F2},~\eqref{eq:F3} is
\begin{equation}
	\phi(t)  =\int_{0}^{t}\mathe^{-\beta (t-\tau)} \chi(\tau) \D\tau,\label{eq:solfilter}
\end{equation}
and
$\upsilon(t) =\int_{0}^{t}\mathe^{-\beta (t-\tau)} u(\tau) \D\tau.$
Integrating the right-hand side of \eqref{eq:solfilter} by parts and using~\eqref{eq:dynIO}, one obtains
that
$\beta \phi(t)  = \chi(t)-\int_{0}^{t}\mathe^{-\beta (t-\tau)}\dot{\chi}(\tau)\D \tau
		 = \chi(t) - \AF \phi(t) - \BF\upsilon(t)-(\int_{0}^{t}\mathe^{-\beta (t-\tau)} \GF \mathe^{\AREF^\top \tau}\D\tau)x_0$.
Under \Cref{ass:noresonance}, since $-\beta\notin \lambda(\AREF)$, the equation $\Gamma \AREF^\top + \beta \Gamma  =\Gamma(\AREF^\top+\beta I_n)= \GF$ admits the unique solution
$\Gamma^\star=\GF(\AREF^\top+\beta I_n)^{-1}$. Hence, letting $\epsilon(t)=\Gamma^\star \mathe^{\AREF^\top t}x_0-\mathe^{-\beta t}\Gamma^\star x_0$,
 one has \[\epsilon(t) =\left( \int_{0}^{t}\mathe^{-\beta (t-\tau)} \GF \mathe^{\AREF^\top \tau}\D \tau\right)x_0.\]
Thus, Eq.~\eqref{eq:deltaAF} holds for $t\geq 0$.

\subsection{Proof of \Cref{lem:null}\label{app:proof5}}

By \cite[Ex.~5.13.16]{meyer2023matrix} and the Cayley-Hamilton theorem, there
exist $\vartheta_{1}^i,\dots,\vartheta_{n}^i\in\R$ such that $\mathe^{\AREF^\top t_i}=\sum_{j=1}^{n} \vartheta_j^i (\AREF^\top)^{j-1}$, for each $i\in\{1,\dots,N\}$.
Therefore, by letting
\begin{equation*}
	\varTheta = \left[\begin{array}{ccccc}
		\vartheta_{1}^{1} -\mathe^{-\beta t_1}& \cdots & \vartheta_{1}^{N}-\mathe^{-\beta t_{N}}\\
		\vdots & \ddots & \vdots\\
		\vartheta_{n}^{1} & \cdots & \vartheta_{n}^{N}
	\end{array}\right]\in\R^{n\times N},
\end{equation*} 
\begin{multline*}
	[\begin{array}{ccc}
		\mathe^{\AREF^\top t_1}-\mathe^{-\beta t_1} I_{n} & \cdots & \mathe^{\AREF^\top t_{N}}-\mathe^{-\beta t_{N}} I_{n}
	\end{array}]\\
	= [\begin{array}{cccc}
		 I_n & \AREF^\top & \cdots & (\AREF^\top)^{n-1}
	\end{array}](\varTheta \otimes I_n).
\end{multline*}
Since $\rank(\varTheta)\leq n$, by the rank plus nullity theorem \cite[Eq.~(4.4.15)]{meyer2023matrix},
$\dim(\Ker(\varTheta))\geq N-n$.
By \cite[Ex.~5.8.15]{meyer2023matrix},  for each $\alpha\in\Ker(\varTheta)$,
one has $(\varTheta \otimes I_n)(\alpha\otimes I_n)=(\varTheta \alpha)\otimes I_n=0_{n^2,n}$.
Therefore, it holds that $	[\begin{array}{ccc}
	\mathe^{\AREF^\top t_1}-\mathe^{-\beta t_1} I_{n} & \cdots & \mathe^{\AREF^\top t_{N}}-\mathe^{-\beta t_{N}} I_{n}
\end{array}](\alpha\otimes I_n)=\sum_{i=1}^{N}\alpha_{i} (\mathe^{\AREF^\top t_i}-\mathe^{-\beta t_i} I_{n})=0_{n,n}$,
for all $\alpha\in\Ker(\varTheta)$.
\subsection{Proof of \Cref{thm:pe}\label{app:proof6}}
Following \cite[App.~A]{11045686}, consider the interconnection of system~\eqref{eq:systemIO} and the filters~\eqref{eq:MIMOfilters},
whose dynamics are 
\begin{align}
	\dot{\xi}(t) &= \AEX \xi(t) + \BEX u(t) + \GEX\wp(t), & \dot{\wp}(t)=\AREF^\top \wp(t)
	\label{eq:series}
\end{align} 
where $\xi(t)=[\begin{array}{ccc}
	\upsilon^\top(t) & \phi^\top(t) & \chi^\top(t)
\end{array}]^\top$ and
\begin{gather*}
\AEX  = \left[\begin{array}{ccc}
	-\beta & 0_{1,2n} & 0_{1,2n} \\
	0_{2n,1} &-\beta I_{2n} & I_{2n}\\
	0_{2n,1} & 0_{2n,2n} & \AF
\end{array}\right], \qquad
\BEX = \left[\begin{array}{c}
	1\\
	0_{2n,1}\\
	\BF
\end{array}\right],\\ 
\GEX = \left[\begin{array}{c}
0_{1,n}\\
0_{2n,n}\\
\GF
\end{array}\right],\qquad
\CEX = \left[\begin{array}{ccc}
0_{2n,1} & I_{2n} & 0_{2n,2n}\\
1 & 0_{1,2n} & 0_{1,2n}
\end{array}\right],
\end{gather*}
whose output $\varkappa(t)=\CEX \xi(t)$ from the initial condition $\xi(0)=0_{4n+1,1}$, $\wp(0)=x_0$,  is
$\varkappa(t)=[\begin{array}{cc}
\phi^\top(t) & \upsilon^\top(t)
\end{array}]^\top$. Under \Cref{ass:noresonance2}, the Sylvester equation
$\Pi \AREF^\top-A^\top \Pi =\BREF C$
admits an unique solution $\Pi^\star\in\R^{n\times n}$. Therefore, the matrix $\PIF^\star=\col(
	0_{n,n} , \Pi^\star
)$ solves 
$\PIF \AREF^\top-\AF \PIF =\GF$ and
hence, 
under \Cref{ass:noresonance},
the matrix $\PIEX^\star=\col(
	0_{1,n},  \PIF^{\star} (\AREF^\top+ \beta I_n)^{- 1} , \PIF^\star)$ solves 
\begin{equation}\label{eq:extendedSylvester}
	\PIEX \AREF^\top-\AEX \PIEX =\GEX.
\end{equation}

By linearity, the solution to system~\eqref{eq:series}  from the initial condition $\xi(0)=0_{4n+1,1}$ satisfies
$\xi(t)=\xi^{\rma}(t)+\xi^{\rmb}(t)$, with
\begin{subequations}
\begin{align}
	\dot{\xi}^{\rma}(t) & = \AEX \xi^{\rma}(t) + \BEX u(t), & \xi^{\rma}(0)&=-\PIEX^\star x_0,\label{eq:undisturbed}\\
	 \dot{\xi}^{\rmb}(t)  & = \AEX \xi^{\rmb}(t) + \GEX\wp(t), & \xi^{\rmb}(0)&=\PIEX^\star x_0.\label{eq:disturbed}
\end{align}
\label{eq:decompositionIntercon}%
\end{subequations}
Since $\wp(t)=\mathe^{\AREF^\top t}x_0$ and $\PIEX^\star$ solves~\eqref{eq:extendedSylvester}, the solution to system~\eqref{eq:disturbed}
is ${\xi}^{\rmb}(t) = \PIEX^\star \mathe^{\AREF^\top t}x_0$.
Therefore, one has \[\varkappa(t)=\CEX {\xi}^{\rma}(t) +\CEX {\xi}^{\rmb}(t)= \varkappa^{\rma}(t) +\CEX  \PIEX^\star \mathe^{\AREF^\top t}x_0.\]
By \cite[App.~A]{11045686}, the dynamics of $ \varkappa^{\rma}(t) $ are given by
\begin{equation}\label{eq:chia}
	\dot{\varkappa}^{\rma}(t) = \AV \varkappa^{\rma}(t)  + \BV u(t) + \GV \eth(t), \; \dot{\eth}(t)=-\beta \eth(t),
\end{equation}
with initial conditions $\varkappa^{\rma}(0)=\varkappa^{\rma}_0$, $\varkappa^{\rma}_0=-\CEX \PIEX^\star x_0$, $\eth(0)=\eth_0$, $\eth_0=\digamma \PIEX^\star x_0$, 
$\digamma=[\begin{array}{ccc}
	\BF & \AF+\beta I_{2n} & -I_{2n}
\end{array}]$,
	$\AV = \left[\begin{smallmatrix}
	\AF & \BF\\
	0_{1, 2n} & - \beta
\end{smallmatrix}\right]$, 
$\BV  = \left[\begin{smallmatrix}
	0_{2n,1} \\ 1
\end{smallmatrix}\right]$, 
$\GV  = \left[\begin{smallmatrix}
	I_{2n}\\
	0_{1,2n}
\end{smallmatrix}\right]$,
and $(\AV,\BV)$ is reachable. Since $(\AF,\BF)$ is reachable, by the PBH test for reachability \cite[Thm.~6.2.6]{kailath1980linear}, $\rank ([\begin{array}{cc}
	\AF+\beta I_{2n} & \BF
\end{array}])=2n$, $\forall \beta\in\R$. Therefore, one has $[\begin{array}{cc}
\AF+\beta I_{2n} & \BF
\end{array}][\begin{array}{cc}
\AF+\beta I_{2n} & \BF
\end{array}]^\dagger=I_{2n}$. Hence, since
	$\AV =\left[\begin{smallmatrix}
		\AF+\beta I_{2n} & \BF\\
		0_{1, 2n} & 0
	\end{smallmatrix}\right]-\beta I_{2n+1}$,
the matrix
 $\PIV^\star=[\begin{array}{cc}
	\AF+\beta I_{2n} & \BF
\end{array}]^\dagger$ solves
\begin{equation}
	\label{eq:betaSylv}
\AV \PIV+\beta \PIV=\GV.
\end{equation}
 By linearity, the solution to system~\eqref{eq:chia}  from the initial condition $\varkappa^{\rma}(0)=\varkappa^{\rma}_0$ satisfies
$\varkappa^{\rma}(t)=\varkappa^{\rma}_{\rmI}(t)+\varkappa^{\rma}_{\rmII}(t)$,
\begin{subequations}
	\begin{align}
		\dot{\varkappa}^{\rma}_{\rmI}(t) & = \AV \varkappa^{\rma}_{\rmI}(t) + \BV u(t), & \varkappa^{\rma}_{\rmI}(0)&=\varkappa^{\rma}_0-\PIV^\star \eth_0,\label{eq:nobeta}\\
		\dot{\varkappa}^{\rma}_{\rmII}(t)  & = \AV \varkappa^{\rma}_{\rmII}(t) + \GV\eth(t), & \varkappa^{\rma}_{\rmII}(0)&=\PIV^\star \eth_0.\label{eq:betasol}
	\end{align}
	\label{eq:decomp}%
\end{subequations}
Since $\PIV^\star$ solves~\eqref{eq:betaSylv}, the solution to~\eqref{eq:betasol} is $\varkappa^{\rma}_{\rmII}(t) =\PIV^\star\mathe^{-\beta t}\eth_0$.
Therefore, the vector $\varkappa(t)$ satisfies
\begin{equation*}
	\varkappa(t) =  \varkappa^{\rma}_{\rmI}(t) +\CEX  \PIEX^\star \mathe^{\AREF^\top t}x_0+\mathe^{-\beta t}\PIV^\star \digamma \PIEX^\star x_0,
\end{equation*}
where $\varkappa^{\rma}_{\rmI}(t) $ solves~\eqref{eq:nobeta}.

Note that $\lambda(\AV)=\{-\beta\}\cup\lambda(\AF)=\{-\beta\}\cup \lambda(A)\cup\lambda(\AREF)$.
Since $(\AV,\BV)$ is reachable, by \cite[Lem.~1]{lopez2022continuous}, if $\{d_k\}_{k=1}^{N-1}$ is persistently exciting of order $4n+2$, 
 $u(t)=d_k$ for all $t\in[k\TS,(k+1)\TS)$, and \eqref{eq:nonpato} holds, then $\rank(H)=4n+2$, 
\begin{equation*}
	H=\left[
		\begin{array}{c}
		\begin{array}{ccccc}
			\varkappa^{\rma}_{\rmI,1} & \cdots & \varkappa^{\rma}_{\rmI,N-2n-1}
		\end{array}\\
		\han{2n+1}(\{d_k\}_{k=1}^{N-1})
	\end{array}
	\right],
\end{equation*}
where $\varkappa^{\rma}_{\rmI,k}=\varkappa^{\rma}_{\rmI}(k\TS)$, $k\in\N$, $k\geq 1$. In particular, one has $\varkappa^{\rma}_{\rmI,k+1} = \AVD \varkappa^{\rma}_{\rmI,k} + \BVD d_k$,
where $\AVD=\mathe^{\AV\TS}$ and $\BVD=(\int_{0}^{\TS}\mathe^{\AVD\tau}\D \tau)\BV$.
Note that, letting $w_0,w_1,\dots,w_{n+1}$ be defined as in \Cref{ex:uniform}, one has
\[\sum_{i=0}^{n+1}w_{n+1-i}(\CEX  \PIEX^\star \mathe^{\AREF^\top t_{i+j}}x_0+\mathe^{-\beta t_{i+j}}\PIV^\star \digamma \PIEX^\star x_0)=0,\]
$j=1,\dots,N-n-1$.
Thus, one has 
\[\left[\begin{array}{c}
	\bar{\Phi}_N\\
	\bar{\Upsilon}_N
\end{array}\right]=[\begin{array}{ccc}
	\sum_{i=0}^{n+1} w_{i}\varkappa^{\rma}_{\rmI,n+2-i} & \cdots & \sum_{i=0}^{n+1} w_{i}\varkappa^{\rma}_{\rmI,N-i}
\end{array}].\]
Hence, by letting $w_i=0$ for $i\in\Z$, $i\geq n+1$, define
$\Lambda = [\begin{smallmatrix}
		 \sum_{i=0}^{2n+1}w_i (\AVD)^{2n+1-i} & \sum_{i=0}^{2n} w_i (\AVD)^{2n-i}\BVD & \cdots & w_0\BVD
	\end{smallmatrix}]$.
By \cite[Thm.~3.2.1]{chen2012optimal}, if~\eqref{eq:nonpato} holds, then $(\AVD,\BVD)$ is reachable.
Therefore, if~\eqref{eq:nonpato} holds, then $\rank(\Lambda)=2n+1$.
Hence, since
\begin{equation*}
[\begin{array}{ccc}
	\sum_{i=0}^{n+1} w_{i}\varkappa^{\rma}_{\rmI,2n+2-i} & \cdots & \sum_{i=0}^{n+1} w_{i}\varkappa^{\rma}_{\rmI,N-i}
\end{array}] =\Lambda H,
\end{equation*}
by the Frobenius inequality \cite[Ex.~4.5.17]{meyer2023matrix}, one has $\rank(\Lambda H)\geq \rank(\Lambda) + \rank(H) - (4n+2)=2n+1$.
Thus, since $\Lambda H\in \R^{(2n+1)\times (N-2n-1)}$,  $\rank(\Lambda H)=2n+1$.

\subsection{Proof of \Cref{thm:interMIMO}\label{app:proof7}}
The proof follows the same lines of that of \Cref{lem:interconnection} given in \Cref{app:proof2} with minor adaptations to deal with the MIMO case.
It is first proved that, $\forall k\in\N$, $k\geq 1$, one has
\begin{equation}
	\label{eq:AINkBINRRR}
	\AMIN^k\BMIN = \left[\begin{array}{c}
		\Au^k\Bu \\
		\sum_{i=0}^{k-1}\Ay^i \By \CM \AM^{k-i-1}\BM\\
		\AM^k \BM\\
	\end{array}\right].
\end{equation}
By definition of $\AMIN$ and $\BMIN$, Eq.~\eqref{eq:AINkBINRRR} holds for $k=1$.
Assuming that it holds for some $k\in\N$, $k\geq 1$, by computing
$\AMIN \cdot \AMIN^k \BMIN$, one has that~\eqref{eq:AINkBINRRR} holds with $k$ substituted by $k+1$.
Therefore, Eq.~\eqref{eq:AINkBINRRR} holds by induction.

Equations~\eqref{eq:AINkBINRRR} and~\eqref{eq:sum} are now used to find a basis for the image of the reachability matrix 
$\RMIN=[\begin{array}{ccc}
	\BMIN & \cdots & \AMIN^{(m+2p)n-1}\BMIN
\end{array}]$ of the pair $(\AMIN,\BMIN)$.
Since $\Au=\AREF\otimes I_m$, $\Bu=\BREF\otimes I_m$, $\Ay=\AREF\otimes I_p$, and $\By=\AREF\otimes I_p$,  by letting $c_0=1$, by~\eqref{eq:AINkBINRRR} and~\eqref{eq:sum}, one has
$\sum_{j=0}^k c_j \AMIN^{k-j}\BMIN = \col(\bbe{n-k}{n}\otimes I_m,\allowbreak \sum_{i=0}^{k-1}(\bbe{n-k+i+1}{n}\otimes I_p)\CM\AM^i\BM, \allowbreak\sum_{j=0}^k c_j \AM^{k-j} \BM)$,~$k=1,\dots,n-1$, and 
$\sum_{j=0}^n c_j \AMIN^{k-j}\BMIN = \col(0_{n,1}\otimes I_m,\allowbreak \sum_{i=0}^{n-1}(\bbe{n-i}{n}\otimes I_p)\CM\AM^{k-i-1}\BM, \allowbreak \sum_{j=0}^n c_j \AM^{k-j} \BM)$,~$k=n,\dots,3n-1$.
Therefore, by letting $\HF\in\R^{(m+2p)n\times (m+2p)n}$ be the upper triangular Toeplitz matrix 
whose first row is $[\begin{array}{cccccccc}
	1 & c_1 & \cdots & c_n & 0_{1,(m+2p-1)n-1}
\end{array}]$, one has
\begin{equation*}
	\RIN (\HF \otimes I_m)= \left[\begin{array}{cc}
		\swap_n\otimes I_m & 0_{nm, (m+2p-1)nm}\\
		M_{\mathrm{M}}( \swap_n\otimes I_m) & O_{\mathrm{M}} R_{\mathrm{M}}\\
		F_{\mathrm{M}}(\swap_n \otimes I_m)& q_{\mathrm{r}}(\AM)R_{\mathrm{M}}
	\end{array}\right],
\end{equation*}
where $q_{\mathrm{r}}(s)=s^n+c_1s^{n-1}+\cdots +c_n$, and
\begin{align*}
	R_{\mathrm{M}} & =[\begin{array}{ccccc}
	\BM &  \cdots & \AM^{(m+2p-1)n}\BM 
\end{array}],\\
F_{\mathrm{M}}& =[\begin{array}{ccc}
	( \sum_{i = 0}^{n - 1} c_i \AM^{n  - i- 1} ) \BM & \cdots & \BM
\end{array}].
\end{align*}
Under \Cref{ass:MIMOass}, the pair $(\AM,\BM)$ is reachable and hence $\rank(R_{\mathrm{M}})=pn$.
Since $O_{\mathrm{M}}$ is invertible and $\swap_n^{-1}=\swap_n$, by taking linear combinations of the columns of $\RIN \HF $, one has 
\begin{equation*}
	\SP(\RMIN)=\SP\left(\left[\begin{array}{cc}
		I_{mn} & 0_{mn,pn}\\
		0_{pn,mn} & I_{pn}\\
		\mho_1 & \mho_2
	\end{array}\right]\right),
\end{equation*} 
where $\mho_1=F_{\mathrm{M}}-q_{\mathrm{r}}(\AM)O_{\mathrm{M}}^{-1}M_{\mathrm{M}}$ and $\mho_2=q_{\mathrm{r}}(\AM)O_{\mathrm{M}}^{-1}$.
Since $(\AM,\BM,\CM)$ is a realization of~\eqref{eq:diffEqMIMO}, one has $\CM\AM^nO_{\mathrm{M}}^{-1}=[\begin{array}{ccc}
	-A_n & \cdots & -A_1
\end{array}]$ and $[\begin{array}{ccc}
	\CM \AM^{n-1}\BM & \cdots & \CM\BM
\end{array}]-\CM\AM^{n}O_{\mathrm{M}}^{-1}M_{\mathrm{M}}=[\begin{array}{ccc}
	B_n & \cdots & B_1
\end{array}]$; see \Cref{rem:comments}. Hence, it results that
$\CM\mho_2 =  [\begin{array}{ccc}
	c_n I_p-A_n & \cdots & c_1I_p-A_1
\end{array}]$ and $\CM\mho_1=[\begin{array}{cccc}
	B_n & \cdots & B_1
\end{array}]$. Finally, since $O_{\mathrm{M}}^{-1}\By=\bbe{1}{n}\otimes I_p$,
one has $\mho_2 \By=\col(
c_nI_p-A_n , \cdots , c_1I_p-A_1
)$.
Thus, let 
\[P=\left[\begin{array}{ccc}
	I_{mn} & 0_{mn,pn} & 0_{mn,pn}\\
	0_{pn,mn} & I_{pn} & 0_{pn,pn}\\
	\mho_1 & \mho_2 & I_{pn}
\end{array}\right].\]
By classical results on the Kalman decomposition for reachability \cite[Sec.~16.2]{hespanha2018linear}, one has
\begin{align*}
	P^{-1}\AMIN P &= \left[\begin{array}{cc}
		\AMINrr & \AMINru\\
		0_{pn,(m+p)n} & \AMINuu
	\end{array}\right], \\ P^{-1}\BMIN & = \left[\begin{array}{c}
		\BMINr \\ 0_{pn,m}
	\end{array}\right],
\end{align*}
and the pair $(\AMINrr,\BMINr)$ is reachable. Note that
\begin{gather*}
	\AMINrr  = \left[\begin{array}{cc}
		\Au & 0_{mn,pn}\\
		\By \CM \mho_1 & \Ay+\By \CM \mho_2
	\end{array}\right],\\
	\BMINr  = \left[\begin{array}{c}
		\Bu \\ 0_{pn,m}
	\end{array}\right],\quad
	\AMINru = \left[\begin{array}{c}
		0_{mn,pn}\\ \By \CM
	\end{array}\right],\\ 
	\AMINuu =\AM- \mho_2 \By \CM=\Ay^\top.
\end{gather*}
Hence, the statement follows by the fact that $\AMINrr  =\AEF$, $\BMINr=\BEF$, $\AINru=\GEF$,
and $\kappa(t)=x(t)-\mho_1\zeta(t)-\mho_2\mu(t)$ satisfies $\dot{\kappa}(t)=\AMINuu \kappa(t)$
with $\kappa(0)=x_0$.
\subsection{Proof of \Cref{lem:relationMIMO}\label{app:proof8}}

Following the same construction employed in Appendix~\ref{app:proof4}, one has that the following relation holds for all $t\geq 0$:
\begin{align*}
	\beta \phi(t)  &= \chi(t)-\int_{0}^{t}\mathe^{-\beta (t-\tau)}\dot{\chi}(\tau)\D \tau\\
		&= \chi(t) - \AEF \phi(t) - \BEF\upsilon(t)\\
		& \qquad -\left(\int_{0}^{t}\mathe^{-\beta (t-\tau)} \GEF \mathe^{\Ay^\top \tau}\D\tau\right)x_0.
\end{align*}
Under \Cref{ass:noresonance}, since $-\beta\notin \lambda(\Ay)$, the equation $\Gamma \Ay^\top + \beta \Gamma  =\Gamma(\Ay^\top+\beta I_{pn})= \GEF$ admits the unique solution
$\Gamma^\star=\GEF(\Ay^\top+\beta I_{pn})^{-1}$. Thus, the statement follows by the fact that 
$(\int_{0}^{t}\mathe^{-\beta (t-\tau)} \GEF \mathe^{\Ay^\top \tau}\D\tau)x_0=\epse(t)$.
\subsection{Proof of \Cref{thm:peMIMO}\label{app:proof9}}
By adapting the construction given in Appendix~\ref{app:proof6} for the SISO case, consider the interconnection of system~\eqref{eq:diffEqMIMO} and the filters~\eqref{eq:MIMOfilters},
whose dynamics are 
\begin{subequations}
\begin{align}
	\dot{\xi}(t) &= \AMEX \xi(t) + \BMEX u(t) + \GMEX\wp(t), \\
	 \dot{\wp}(t)&=\Ay^\top \wp(t)
\end{align} 
\label{eq:seriesMIMO}%
\end{subequations}
where $\xi(t)=[\begin{array}{ccc}
	\upsilon^\top(t) & \phi^\top(t) & \chi^\top(t)
\end{array}]^\top$ and
\begin{gather*}
	\AMEX  = \left[\begin{array}{ccc}
		-\beta I_m& 0_{m,(m+p)n} & 0_{m,(m+p)n} \\
		0_{(m+p)n,m} &-\beta I_{(m+p)n} & I_{(m+p)n}\\
		0_{(m+p)n,m} & 0_{(m+p)n,(m+p)n} & \AEF
	\end{array}\right], \\
	\BMEX = \left[\begin{array}{c}
		I_m\\
		0_{(m+p)n,m}\\
		\BEF
	\end{array}\right],\\ 
	\GMEX = \left[\begin{array}{c}
		0_{pn}\\
		0_{(m+p)n,pn}\\
		\GEF
	\end{array}\right],\\
	\CMEX = \left[\begin{array}{ccc}
		0_{(m+p)n,m} & I_{(m+p)n} & 0_{(m+p)n,(m+p)n}\\
		I_m & 0_{m,(m+p)n} & 0_{m,(m+p)n}
	\end{array}\right],
\end{gather*}
whose output $\varkappa(t)=\CMEX \xi(t)$ from the initial condition $\xi(0)=0_{2(m+p)n+m,1}$, $\wp(0)=x_0$,  is
$\varkappa(t)=[\begin{array}{cc}
	\phi^\top(t) & \upsilon^\top(t)
\end{array}]^\top$. 
If $\sigma(\AM)\cap \sigma(\AREF)=\emptyset$,
then the Sylvester equation
$\PiM \Ay^\top-\AM^\top \PiM =\By \CM$
admits an unique solution $\PiM^\star\in\R^{pn\times pn}$. Therefore, the matrix $\PIEF^\star=\col(
0_{mn,pn} , \PiM^\star
)$ solves 
$\PIEF \Ay^\top-\AEF \PIEF =\GEF$ and
hence, 
under \Cref{ass:noresonance},
the matrix $\PIMEX^\star=\col(
0_{m,pn},  \PIEF^{\star} (\Ay^\top+ \beta I_{pn})^{- 1} , \PIEF^\star)$ solves 
\begin{equation}\label{eq:extendedSylvesterMIMO}
	\PIMEX \Ay^\top-\AMEX \PIMEX =\GMEX.
\end{equation}

By linearity, the solution to~\eqref{eq:seriesMIMO}  from the initial condition $\xi(0)=0_{2(m+p)n+m,1}$ satisfies
$\xi(t)=\xi^{\rma}(t)+\xi^{\rmb}(t)$, with
\begin{subequations}
	\begin{align}
		\dot{\xi}^{\rma}(t) & = \AMEX \xi^{\rma}(t) + \BMEX u(t), & \xi^{\rma}(0)&=-\PIMEX^\star x_0,\label{eq:undisturbedMIMO}\\
		\dot{\xi}^{\rmb}(t)  & = \AMEX \xi^{\rmb}(t) + \GMEX\wp(t), & \xi^{\rmb}(0)&=\PIMEX^\star x_0.\label{eq:disturbedMIMO}
	\end{align}
	\label{eq:decompMIMOositionInterconMIMO}%
\end{subequations}
Since $\wp(t)=\mathe^{\Ay^\top t}x_0$ and $\PIMEX^\star$ solves~\eqref{eq:extendedSylvesterMIMO}, the solution to system~\eqref{eq:disturbedMIMO}
is ${\xi}^{\rmb}(t) = \PIMEX^\star \mathe^{\Ay^\top t}x_0$.
Therefore, one has \[\varkappa(t)=\CMEX {\xi}^{\rma}(t) +\CMEX {\xi}^{\rmb}(t)= \varkappa^{\rma}(t) +\CMEX  \PIMEX^\star \mathe^{\Ay^\top t}x_0.\]
By \cite[App.~A]{11045686}, the dynamics of $ \varkappa^{\rma}(t) $ are given by
\begin{subequations}
\begin{align}
	\dot{\varkappa}^{\rma}(t) &= \AMV \varkappa^{\rma}(t)  + \BMV u(t) + \GMV \eth(t),\\ \dot{\eth}(t)&=-\beta \eth(t),
\end{align}
\label{eq:chiaMIMO}%
\end{subequations}
with initial conditions $\varkappa^{\rma}(0)=\varkappa^{\rma}_0$, $\varkappa^{\rma}_0=-\CMEX \PIMEX^\star x_0$, $\eth(0)=\eth_0$, $\eth_0=\digamma_{\mathrm{M}} \PIMEX^\star x_0$, 
$\digamma_{\mathrm{M}}=[\begin{array}{ccc}
	\BEF & \AEF+\beta I_{(m+p)n} & -I_{(m+p)n}
\end{array}]$,
$\AMV = \left[\begin{smallmatrix}
	\AEF & \BEF\\
	0_{m, (m+p)n} & - \beta I_m
\end{smallmatrix}\right]$, 
$\BMV  = \left[\begin{smallmatrix}
	0_{(m+p)n,m} \\ I_m
\end{smallmatrix}\right]$, 
$\GMV  = \left[\begin{smallmatrix}
	I_{(m+p)n}\\
	0_{m,(m+p)n}
\end{smallmatrix}\right]$,
and $(\AMV,\BMV)$ is reachable. Since $(\AEF,\BEF)$ is reachable, by the PBH test for reachability \cite[Thm.~6.2.6]{kailath1980linear}, $\rank ([\begin{array}{cc}
	\AEF+\beta I_{(m+p)n} & \BEF
\end{array}])=(m+p)n$, $\forall \beta\in\R$. Therefore, one has $[\begin{array}{cc}
	\AEF+\beta I_{(m+p)n} & \BEF
\end{array}][\begin{array}{cc}
	\AEF+\beta I_{(m+p)n} & \BEF
\end{array}]^\dagger=I_{(m+p)n}$. Hence, since
$\AMV =\left[\begin{smallmatrix}
	\AEF+\beta I_{(m+p)n} & \BEF\\
	0_{m, (m+p)n} & 0_{m,m}
\end{smallmatrix}\right]-\beta I_{(m+p)n+m}$,
$\PIMV^\star=[\begin{array}{cc}
	\AEF+\beta I_{(m+p)n} & \BEF
\end{array}]^\dagger$ solves
\begin{equation}
	\label{eq:betaSylvMIMO}
	\AMV \PIMV+\beta \PIMV=\GMV.
\end{equation}
By linearity, the solution to system~\eqref{eq:chiaMIMO}  from the initial condition $\varkappa^{\rma}(0)=\varkappa^{\rma}_0$ satisfies
$\varkappa^{\rma}(t)=\varkappa^{\rma}_{\rmI}(t)+\varkappa^{\rma}_{\rmII}(t)$,
\begin{subequations}
	\begin{align}
		\dot{\varkappa}^{\rma}_{\rmI}(t) &= \AMV \varkappa^{\rma}_{\rmI}(t) + \BMV u(t),  \label{eq:nobetaMIMO}\\
		\dot{\varkappa}^{\rma}_{\rmII}(t)   &= \AMV \varkappa^{\rma}_{\rmII}(t) + \GMV\eth(t),\label{eq:betasolMIMO}
	\end{align}
	\label{eq:decompMIMO}%
\end{subequations}
where $\varkappa^{\rma}_{\rmI}(0)=\varkappa^{\rma}_0-\PIMV^\star \eth_0$ and $\varkappa^{\rma}_{\rmII}(0)=\PIMV^\star \eth_0$.
Since $\PIMV^\star$ solves~\eqref{eq:betaSylvMIMO}, the solution to~\eqref{eq:betasolMIMO} is $\varkappa^{\rma}_{\rmII}(t) =\PIMV^\star\mathe^{-\beta t}\eth_0$.
Therefore, the vector $\varkappa(t)$ satisfies
\begin{equation*}
	\varkappa(t) =  \varkappa^{\rma}_{\rmI}(t) +\CMEX  \PIMEX^\star \mathe^{\Ay^\top t}x_0+\mathe^{-\beta t}\PIMV^\star \digamma_{\mathrm{M}} \PIMEX^\star x_0,
\end{equation*}
where $\varkappa^{\rma}_{\rmI}(t) $ solves~\eqref{eq:nobetaMIMO}.

Note that $\lambda(\AMV)=\{-\beta\}\cup\lambda(\AEF)=\{-\beta\}\cup \lambda(\AM)\cup\lambda(\AREF)$.
Since $(\AMV,\BMV)$ is reachable, by \cite[Lem.~1]{lopez2022continuous}, if $\{d_k\}_{k=1}^{N-1}$ is persistently exciting of order $2((m+p)n+m)$, 
$u(t)=d_k$ for all $t\in[k\TS,(k+1)\TS)$, and \eqref{eq:nonpatoMIMO} holds, then $\rank(H)=(m+p)n+m + m((m+p)n+m)$, 
\begin{equation*}
	H=\left[
	\begin{array}{c}
		\begin{array}{ccccc}
			\varkappa^{\rma}_{\rmI,1} & \cdots & \varkappa^{\rma}_{\rmI,N-((m+p)n+m)}
		\end{array}\\
		\han{(m+p)n+m}(\{d_k\}_{k=1}^{N-1})
	\end{array}
	\right],
\end{equation*}
where $\varkappa^{\rma}_{\rmI,k}=\varkappa^{\rma}_{\rmI}(k\TS)$, $k\in\N$, $k\geq 1$. In particular, one has $\varkappa^{\rma}_{\rmI,k+1} = \AMVD \varkappa^{\rma}_{\rmI,k} + \BMVD d_k$,
where $\AMVD=\mathe^{\AMV\TS}$ and $\BMVD=(\int_{0}^{\TS}\mathe^{\AMVD\tau}\D \tau)\BMV$.
Note that, letting $w_0,w_1,\dots,w_{n+1}$ be defined as in \Cref{ex:uniform}, one has
$\sum_{i=0}^{n+1}w_{n+1-i}(\CMEX  \PIMEX^\star \mathe^{\Ay^\top t_{i+j}}x_0+\mathe^{-\beta t_{i+j}}\PIMV^\star \digamma_{\mathrm{M}} \PIMEX^\star x_0)=0$,
$j=1,\dots,N-n-1$.
Thus, one has that
\[\left[\begin{array}{c}
	\bar{\Phi}_N\\
	\bar{\Upsilon}_N
\end{array}\right]=[\begin{array}{ccc}
	\sum_{i=0}^{n+1} w_{i}\varkappa^{\rma}_{\rmI,n+2-i} & \cdots & \sum_{i=0}^{n+1} w_{i}\varkappa^{\rma}_{\rmI,N-i}
\end{array}].\]
Hence, by letting $w_i=0$ for $i\in\Z$, $i\geq n+1$, define $\Lambda$ as in Eq.~\eqref{eq:Lambda}.
\begin{figure*}[b!]
	\hrule
	\begin{equation}\label{eq:Lambda}
		\Lambda = \left[\begin{array}{cccc}
			\sum_{i=0}^{(m+p)n+m}w_i (\AMVD)^{(m+p)n+m-i}  & \sum_{i=0}^{(m+p)n+m-1} w_i (\AMVD)^{(m+p)n+m-1-i}\BMVD  & \cdots & w_0\BMVD
		\end{array}
		\right].
	\end{equation}
\end{figure*}
By \cite[Thm.~3.2.1]{chen2012optimal}, if~\eqref{eq:nonpatoMIMO} holds, then $(\AMVD,\BMVD)$ is reachable.
Therefore, if~\eqref{eq:nonpatoMIMO} holds, then $\rank(\Lambda)=(m+p)n+m$.
Hence, since
\begin{multline*}
	[\begin{array}{ccc}
		\sum_{i=0}^{n+1} w_{i}\varkappa^{\rma}_{\rmI,(m+p)n+m+1-i} & \cdots & \sum_{i=0}^{n+1} w_{i}\varkappa^{\rma}_{\rmI,N-i}
	\end{array}] \\=\Lambda H,
\end{multline*}
by the Frobenius inequality \cite[Ex.~4.5.17]{meyer2023matrix}, one has $\rank(\Lambda H)\geq \rank(\Lambda) + \rank(H) - ((m+p)n+m + m((m+p)n+m))=(m+p)n+m$.
Thus, one has that $\rank(\Lambda H)=(m+p)n+m$.

\bibliographystyle{ieeetr}
\bibliography{biblio}

\end{document}